\documentclass[11pt]{article}
\usepackage{amsfonts,amssymb,amsmath,mathrsfs,bbm}
\usepackage{color,latexsym,amsfonts}

\numberwithin{equation}{section}
\newtheorem{theorem}{Theorem}[section]
\newtheorem{lemma}{Lemma}

\newtheorem{corollary}[theorem]{Corollary}

\def\<{\langle}
\def\>{\rangle}

\begin{document}

\bigskip \bigskip \noindent {\Large \textbf{Reduced critical processes for
small populations }}\footnote{%
\noindent Supported by NSFC (NO.11531001, 11626245), High-end Foreign
Experts Recruitment Program (GDW20171100029) and the Russian Science
Foundation under the grant 14-50-00005}

\noindent {%\normalsize\sf
Minzhi Liu\footnote{%
School of Mathematical Sciences \& Laboratory of Mathematics and Complex
Systems, Beijing Normal University, Beijing 100875, P.R. China. Email:
liuminzhi@mail.bnu.edu.cn}} ~ Vladimir Vatutin \footnote{%
Steklov Mathematical Institute, Gubkin street, 117 996 Moscow, Russia \&
Beijing Normal University, Beijing 100875, P.R. China Email:
vatutin@mi.ras.ru}

\begin{center}
\begin{minipage}{12cm}
\begin{center}\textbf{Abstract}\end{center}
\footnotesize
Let $\left\{ Z(n),n\geq 1\right\} $ be a critical Galton-Watson branching
process with finite variance for the offspring size of particles. Assuming that $0<Z(n)\leq
\varphi (n)$, where either $\varphi (n)=an$ for some $a>0$ or $\varphi
(n)=o(n)$ as $n\rightarrow \infty $, we study the structure of the process $%
\left\{ Z(m,n),0\leq m\leq n\right\} ,$ where $Z(m,n)$ is the number of
particles in the process at moment $m\leq n$ \ having a positive number of
descendants at moment $n$.

\bigskip

\mbox{}\textbf{Keywords:}\quad critical branching process, reduced processes, conditional limit theorem; \\
\mbox{}\textbf{Mathematics Subject Classification}:  Primary 60J80;
secondary 60G50.

\end{minipage}
\end{center}

\section{ Introduction and main results\label{Introd}}

Let $\left\{ Z(n),n\geq 0\right\} $ be a Galton-Watson branching process
with $Z(0)=1$ in which particles produce children in accordance with
probability generating function
\begin{equation*}
f(s)=\mathbf{E}s^{\xi }=\sum_{k=0}^{\infty }f_{k}s^{k}
\end{equation*}%
and let $Z(m,n)$ be the number of particles in the process at moment $m\leq n
$ \ having a positive number of descendants at moment $n$. The process $%
\left\{ Z(m,n),0\leq m\leq n\right\} $ is called a \textit{reduced process}.

Reduced processes for ordinary Galton--Watson branching processes were
introduced by Fleischmann and Prehn \cite{FP}, who discussed the subcritical
case. The distance to the most recent common ancestor (MRCA) for the
supercritical Galton--Watson processes and for the critical processes with
possibly infinite variance of the offspring size has been investigated by
Zubkov \cite{Zub}. Fleischmann and Siegmund-Schultze \cite{FZ} proved a
functional conditional limit theorem establishing, under the condition $%
\left\{ Z(n)>0\right\} $ convergence of the reduced critical Galton--Watson
branching process to the Yule process. \ Different questions related to the
problem of the distribution of the MRCA for the $k$ particles selected at
random among the $Z(n)\geq k$ particles existing in the population at moment
$n$ were considered, for instance, in \cite{Athr2012},\cite{Ath2012b},\cite%
{Dur78},\cite{HJR17}-\cite{Con95}.

However, all these papers do not consider the situation when the size of the
population at moment $n$ is bounded from above. In the present paper, we
study the structure of a critical reduced process and investigate the
asymptotic behavior of the number of its particles under the condition that
the size of the population is bounded and positive at the moment of
observation. Note that the critical Galton-Watson process given its
extinction moment is fixed was investigated in \cite{Sen67} for the
single-type case and in \cite{VD2015} for the multitype setting.

It is known (see, for instance, \cite{AN72}, Chapter I, Section 9 or \cite%
{Sev74}, Chapter II, Section 5) \ that if
\begin{equation}
\mathbf{E}\xi =1,\quad 2B:=Var\xi \in \left( 0,\infty \right) ,
\label{BasicCond}
\end{equation}%
then
\begin{equation}
Q(n):=\mathbf{P}\left( Z(n)>0\right) \sim \frac{1}{Bn}\text{ \ \ as \ \ }%
n\rightarrow \infty   \label{SurvivalProbab}
\end{equation}%
and, for any $y\geq 0$%
\begin{equation}
\lim_{n\rightarrow \infty }\mathbf{P}\left( \frac{Z(n)}{Bn}\leq
y|Z(n)>0\right) =1-e^{-y}.  \label{Yaglom}
\end{equation}%
In addition (see \cite{FZ}), for any fixed $t\in \lbrack 0,1)$ and all $s\in %
\left[ 0,1\right] $
\begin{equation}
\lim_{n\rightarrow \infty }\mathbf{E}\left[ s^{Z(nt,n)}|Z(n)>0\right] =s%
\frac{1-t}{1-ts}.  \label{MRCA_ordinary}
\end{equation}

In this note we study the asymptotic properties of the reduced process when
the condition $\left\{ Z(n)>0\right\} $ is replaced either by the assumption
that $\left\{ 0<Z(n)\leq B\varphi (n)\right\} $ for a function $\varphi
(n)=o(n)$ as $n\rightarrow \infty $ or by the assumption that $\left\{
0<Z(n)\leq aBn\right\} $ for some $a>0$ $.$ Our main results are contained
in two theorems which we formulate below.

\begin{theorem}
\label{T_main}If g.c.d.$\left\{ k:f_{k}>0\right\} =1,$ condition (\ref%
{BasicCond}) is valid, and $\varphi (n)\rightarrow \infty $ in such a way
that $\varphi (n)=o(n)$, then for any $x\in \left( 0,\infty \right) $%
\begin{equation*}
\lim_{n\rightarrow \infty }\mathbf{E}\left[ s^{Z(n-x\varphi (n),n)}\big|%
0<Z(n)\leq B\varphi (n)\right] =sx\frac{1-e^{-(1-s)/x}}{1-s}.
\end{equation*}
\end{theorem}

Let
\begin{equation*}
\beta (n):=\max \left( 0\leq m<n:Z(m,n)=1\right)
\end{equation*}%
be the birth moment of the MRCA of all particles existing in the population
at moment $n$ and let $d(n):=n-\beta (n)$ be the distance from the point of
observation $n$ to the birth moment of the MRCA.

\begin{corollary}
\label{Cor1}Under the conditions of Theorem \ref{T_main}%
\begin{equation*}
\lim_{n\rightarrow \infty }\mathbf{P}\left( d(n)\leq x\varphi (n)|0<Z(n)\leq
B\varphi (n)\right) =x\left( 1-e^{-1/x}\right) .
\end{equation*}
\end{corollary}

To give a complete description of possible situations we present the
following statement.

\begin{theorem}
\label{T_simple}If g.c.d.$\left\{ k:f_{k}>0\right\} =1$ and condition (\ref%
{BasicCond}) is valid, then, for any fixed $t\in \lbrack 0,1)$ and any $a>0$%
\begin{equation*}
\lim_{n\rightarrow \infty }\mathbf{E}\left[ s^{Z(nt,n)}\big|0<Z(n)\leq aBn%
\right] =s\frac{1-t}{1-ts}\frac{1-e^{-\left( 1-ts\right) a/(1-t)}}{1-e^{-a}}.
\end{equation*}
\end{theorem}

\begin{corollary}
\label{Cor2}Under the conditions of Theorem \ref{T_simple}%
\begin{equation*}
\lim_{n\rightarrow \infty }\mathbf{P}\left( d(n)\leq tn|0<Z(n)\leq
aBn\right) =t\frac{1-e^{-a/t}}{1-e^{-a}}.
\end{equation*}
\end{corollary}

\section{Proof of Theorem \protect\ref{T_main}}

For convenience of references we recall Fa\`a di Bruno's formula for the
derivatives of composite functions:

If $i_{r}\in \mathbb{N}_{0}:=\mathbb{N}\cup \left\{ 0\right\} ,r=1,2,...,k$,
$I_{k}:=i_{1}+\cdots +i_{k}$ and
\begin{equation*}
\mathcal{D}(k):=\left\{ \left( i_{1},...,i_{k}\right) :1\cdot i_{1}+2\cdot
i_{2}+\cdot \cdot \cdot +ki_{k}=k\right\} ,
\end{equation*}%
then%
\begin{equation*}
\frac{d^{k}}{d^{k}z}\left[ F(G(z))\right] =\sum_{\mathcal{D}(k)}\frac{k!}{%
i_{1}!\cdot \cdot \cdot i_{k}!}F^{(I_{k})}(G(z))\prod_{r=1}^{k}\left( \frac{%
G^{(r)}(z)}{r!}\right) ^{i_{r}}.
\end{equation*}

We split the proof of Theorem \ref{T_main} into several lemmas.

Let%
\begin{equation*}
f_{0}(s):=s\text{ and }f_{n+1}(s):=f(f_{n}(s)),n\geq 0.
\end{equation*}

Below for arbitrary $x>0$ we agree consider $f_{xn}(s)$ as $f_{\left[ xn%
\right] }(s).$ Besides, the symbol $\sim $ will be usually used (if no
otherwise is stated) for $\overset{n\rightarrow \infty }{\sim }$.

\begin{lemma}
\label{L_Deriv_n}If condition (\ref{BasicCond}) is valid then, for any fixed
$k\in \mathbb{N}:=\left\{ 1,2,...\right\} $ and any fixed $x\in \left(
0,\infty \right) $%
\begin{equation}
f_{n}^{(k)}\left( f_{xn}(0)\right) \sim \frac{k!x^{2}\left( Bxn\right)^{k-1}%
}{\left( x+1\right)^{k+1}}\text{ as }n\rightarrow \infty \text{.}
\label{Der_fn}
\end{equation}
\end{lemma}

\textbf{Proof}. In view of (\ref{SurvivalProbab}) we may rewrite (\ref%
{Yaglom}) in terms of probability functions and Laplace transforms as
follows: for any $\lambda >0$%
\begin{eqnarray}
\lim_{n\rightarrow \infty }\mathbf{E}\left[ e^{-\lambda Z(n)/Bn}|Z(n)>0%
\right] &=&\lim_{n\rightarrow \infty }\frac{f_{n}(e^{-\lambda /Bn})-f_{n}(0)%
}{1-f_{n}(0)}  \notag \\
&=&\lim_{n\rightarrow \infty }Bn\left( f_{n}(e^{-\lambda
Q(n)})-f_{n}(0)\right) =\frac{1}{1+\lambda }.  \label{Laplace0}
\end{eqnarray}
Since
\begin{equation}
\log f_{xn}(0)\sim -(1-f_{xn}(0))\sim -1/Bxn\text{ as }n\rightarrow \infty ,
\label{logAsimpt}
\end{equation}
we conclude that%
\begin{eqnarray*}
\lim_{n\rightarrow \infty }Bn\left( f_{n}(f_{xn}^{\lambda
}(0))-f_{n}(0)\right) &=&\lim_{n\rightarrow \infty }Bn\left( f_{n}\left(
e^{-\lambda (1-f_{xn}(0))}\right) -f_{n}(0)\right) \\
&=&\frac{x}{x+\lambda }.
\end{eqnarray*}%
Clearly, the prelimiting and limiting functions in the previous relations
are analytical in the complex semi-plane Re $\lambda >0.$ Therefore, the
derivatives of any order of the prelimiting functions converge to the
respective derivatives of the limiting function for each $\lambda $ with Re $%
\lambda >0$. Thus,%
\begin{equation*}
\lim_{n\rightarrow \infty }Bn\frac{d^{k}f_{n}(f_{xn}^{\lambda }(0))}{%
d^{k}\lambda }=\left( -1\right) ^{k}\frac{k!x}{\left( x+\lambda \right)
^{k+1}}.
\end{equation*}%
In particular,%
\begin{equation*}
Bn\frac{df_{n}(f_{xn}^{\lambda }(0))}{d\lambda }=Bnf_{n}^{\prime
}(f_{xn}^{\lambda }(0))f_{xn}^{\lambda }(0)\log f_{xn}(0)\sim \left(
-1\right) \frac{1!x}{\left( x+\lambda \right) ^{2}}.
\end{equation*}%
Hence, setting $\lambda =1$ and taking into account (\ref{logAsimpt}) we
conclude that%
\begin{equation*}
f_{n}^{\prime }(f_{xn}(0))\sim \frac{1!x^{2}}{\left( x+1\right) ^{2}},
\end{equation*}%
proving the lemma for $k=1$. Assume that (\ref{Der_fn}) is proved for all $%
k<j$. Since%
\begin{equation*}
\frac{d^{r}}{d^{r}\lambda }f_{xn}^{\lambda }(0)=f_{xn}^{\lambda }(0)\log
^{r}f_{xn}(0),r=1,2,\ldots ,
\end{equation*}%
using Fa\`{a} di Bruno's formula and induction hypothesis we get%
\begin{eqnarray*}
Bn\frac{d^{j}f_{n}(f_{xn}^{\lambda }(0))}{d^{j}\lambda } &=&Bn\sum_{\mathcal{%
D}(j)}\frac{j!}{i_{1}!\cdot \cdot \cdot i_{j}!}f_{n}^{(I_{j})}(f_{xn}^{%
\lambda }(0))\prod_{r=1}^{j}\left( \frac{1}{r!}\frac{d^{r}}{d^{r}\lambda }%
f_{xn}^{\lambda }(0)\right)^{i_{r}} \\
&= &Bn\log^{j}f_{xn}(0)\sum_{\mathcal{D}(j)}\frac{j!}{i_{1}!\cdot \cdot
\cdot i_{j}!}f_{n}^{(I_{j})}(f_{xn}^{\lambda }(0))\prod_{r=1}^{j}\left(
\frac{1}{r!}f_{xn}^{\lambda }(0)\right) ^{i_{r}} \\
&\sim &\left( -1\right)^{j}\frac{Bn}{\left( Bxn\right)^{j}}%
f_{n}^{(j)}(f_{xn}^{\lambda }(0))\sim \left( -1\right)^{j}\frac{j!x}{\left(
x+\lambda \right)^{j+1}}.
\end{eqnarray*}
Hence, setting $\lambda =1$ we obtain%
\begin{equation*}
f_{n}^{(j)}(f_{xn}(0))\sim \frac{j!x^{2}\left( Bxn\right) ^{j-1}}{\left(
x+1\right) ^{j+1.}}
\end{equation*}%
justifying the induction step.

Lemma \ref{L_Deriv_n} is proved. \hfill\rule{2mm}{3mm}\vspace{4mm}

\begin{lemma}
\label{L_derivative}If condition (\ref{BasicCond}) is valid, $m=n-x\varphi
(n)$, where $x\in \left( 0,\infty \right) $ and $\varphi (n)=o(n)$ as $%
n\rightarrow \infty $, then, for any fixed $j\in \mathbb{N}:=\left\{
1,2,...\right\} $
\begin{equation*}
f_{m}^{(j)}\left( f_{x\varphi (n)}(0)\right) \sim \frac{j!\left( Bx\varphi
(n)\right)^{j+1}}{B^{2}n^{2}}\text{ as }n\rightarrow \infty \text{.}
\end{equation*}
\end{lemma}

\textbf{Proof}. It is known (see, for instance \cite{AN72}, Chapter 1,
Section 9, Corollary 1) that under condition (\ref{BasicCond})
\begin{equation*}
\lim_{n\rightarrow \infty }n^{2}\left[ f_{n+1}(0)-f_{n}(0)\right] =\frac{1}{B%
}.
\end{equation*}

We consider for $\lambda >0$ the function%
\begin{equation*}
f_{m}(f_{x\varphi (n)}^{\lambda }(0))=f_{m}(e^{\lambda \log f_{x\varphi
(n)}(0)})
\end{equation*}%
and find $r$ such that%
\begin{equation*}
1-f_{r+1}(0)<1-f_{x\varphi (n)}^{\lambda }(0)\leq 1-f_{r}(0).
\end{equation*}%
In view of (\ref{SurvivalProbab}) we know that
\begin{equation*}
1-f_{x\varphi (n)}^{\lambda }(0)\sim \lambda \left( 1-f_{x\varphi
(n)}(0)\right) \sim \frac{\lambda }{Bx\varphi (n)}.
\end{equation*}%
Hence we get%
\begin{equation*}
r\sim \frac{x\varphi (n)}{\lambda }=o(n)\text{ as }n\rightarrow \infty \text{%
.}
\end{equation*}%
Then for $n-m=x\varphi (n)$%
\begin{eqnarray*}
&&\lim_{n\rightarrow \infty }\frac{n^{2}}{x\varphi (n)}\left[
f_{m}(f_{r}(0))-f_{m}(0)\right] \\
&=&\lim_{n\rightarrow \infty }\frac{1}{x\varphi (n)}\sum_{k=0}^{r-1}n^{2}%
\left[ f_{m}(f_{k+1}(0))-f_{m}(f_{k}(0))\right] \\
&=&\lim_{n\rightarrow \infty }\frac{1}{x\varphi (n)}\sum_{k=0}^{r-1}\frac{%
n^{2}}{\left( m+k\right)^{2}}\left( m+k\right)^{2}\left[
f_{m+k+1}(0)-f_{m+k}(0)\right] \\
&=&\frac{1}{B}\lim_{n\rightarrow \infty }\frac{1}{x\varphi (n)}%
\sum_{k=0}^{r-1}1=\frac{1}{B}\frac{1}{\lambda }.
\end{eqnarray*}
Thus,%
\begin{equation}
\lim_{n\rightarrow \infty }\frac{n^{2}}{x\varphi (n)}\left[ f_{m}(e^{\lambda
\log f_{x\varphi (n)}(0)})-f_{m}(0)\right] =\frac{1}{B}\frac{1}{\lambda }%
,\quad \lambda >0.  \label{Complex}
\end{equation}%
According to the similar reason in the proof of Lemma \ref{L_Deriv_n}, we
have that for each $k\geq 1$%
\begin{equation}
\lim_{n\rightarrow \infty }\frac{Bn^{2}}{x\varphi (n)}\frac{d^{k}}{%
d^{k}\lambda }\left[ f_{m}(e^{\lambda \log f_{x\varphi (n)}(0)})\right]
=\left( -1\right) ^{k}\frac{k!}{\lambda^{k+1}}.  \label{Derivative2}
\end{equation}

By Fa\`a di Bruno's formula we have%
\begin{eqnarray*}
&&\frac{d^{k}}{d^{k}\lambda }\left[ f_{m}(e^{\lambda \log f_{x\varphi
(n)}(0)})\right] \\
&=&\sum_{\mathcal{D}(k)}\frac{k!}{i_{1}!\cdots i_{k}!}f_{m}^{(I_{k})}(e^{%
\lambda \log f_{x\varphi (n)}(0)})\prod_{r=1}^{k}\left( \left( \frac{%
e^{\lambda \log f_{x\varphi (n)}(0)}}{r!}\right)^{(r)}\right)^{i_{r}} \\
&=&\sum_{\mathcal{D}(k)}\frac{k!}{i_{1}!\cdots i_{k}!}f_{m}^{(I_{k})}(e^{%
\lambda \log f_{x\varphi (n)}(0)})e^{\lambda I_{k}\log f_{x\varphi
(n)}(0)}\prod_{r=1}^{k}\frac{\left( \log f_{x\varphi (n)}(0)\right) }{\left(
r!\right)^{i_{r}}}^{ri_{r}} \\
&=&\left( \log f_{x\varphi (n)}(0)\right)^{k}\sum_{\mathcal{D}(k)}\frac{k!}{%
i_{1}!\cdots i_{k}!}f_{m}^{(I_{k})}(e^{\lambda \log f_{x\varphi
(n)}(0)})e^{\lambda I_{k}\log f_{x\varphi (n)}(0)}\prod_{r=1}^{k}\left(
\frac{1}{r!}\right)^{i_{r}}.
\end{eqnarray*}%
Recalling (\ref{logAsimpt}) we get
\begin{eqnarray*}
&&\frac{Bn^{2}}{x\varphi (n)}\frac{d^{k}}{d^{k}\lambda }\left[
f_{m}(e^{\lambda \log f_{x\varphi (n)}(0)})\right] \left\vert _{\lambda
=1}\right. \\
&\sim &\left( -1\right)^{k}\sum_{\mathcal{D}(k)}\frac{k!}{i_{1}!\cdots i_{k}!%
}B^{2}n^{2}\frac{f_{m}^{(I_{k})}(f_{x\varphi (n)}(0))}{(Bx\varphi (n))^{k+1}}%
\prod_{r=1}^{k}\left( \frac{1}{r!}\right)^{i_{r}}\sim \left( -1\right)^{k}k!.
\end{eqnarray*}

In particular,%
\begin{equation*}
\frac{Bn^{2}}{x\varphi (n)}\frac{d}{d\lambda }\left[ f_{m}(e^{\lambda \log
f_{x\varphi (n)}(0)})\right] \left\vert _{\lambda =1}\right. \sim -\frac{%
n^{2}}{\left( x\varphi (n)\right) ^{2}}f_{m}^{\prime }(f_{x\varphi
(n)}(0))\sim \left( -1\right) 1!
\end{equation*}%
giving
\begin{equation*}
f_{m}^{\prime }(f_{x\varphi (n)}(0))\sim \frac{\left( x\varphi (n)\right)
^{2}}{n^{2}}.
\end{equation*}%
Now, by induction we prove that, for any $k\geq 1$, as $n\rightarrow \infty $%
,
\begin{equation*}
\frac{B^{2}n^{2}}{\left( xB\varphi (n)\right) ^{k+1}}f_{m}^{(k)}\left(
f_{x\varphi (n)}(0)\right) \sim k!.
\end{equation*}%
This is true for $k=1$ and if this is true for $k<j$ then, in view of (\ref%
{Derivative2}) and the induction hypothesis
\begin{eqnarray*}
&&\frac{Bn^{2}}{x\varphi (n)}\frac{d^{j}}{d^{j}\lambda }\left[
f_{m}(e^{\lambda \log f_{x\varphi (n)}(0)})\right] \left\vert _{\lambda
=1}\right. \\
&\sim &\left( -1\right) ^{j}\sum_{\mathcal{D}(j)}\frac{j!}{i_{1}!i_{2}!\cdot
\cdot \cdot i_{j}!}\frac{B^{2}n^{2}}{(Bx\varphi (n))^{j+1}}%
f_{m}^{(I_{j})}(f_{x\varphi (n)}(0))\prod_{r=1}^{j}\frac{1}{\left( r!\right)
^{i_{r}}} \\
&\sim &\left( -1\right) ^{j}\frac{j!}{j!0!\cdot \cdot \cdot 0!}\frac{%
B^{2}n^{2}}{(Bx\varphi (n))^{j+1}}f_{m}^{(j)}(f_{x\varphi (n)}(0)) \\
&=&\left( -1\right) ^{j}\frac{B^{2}n^{2}}{(Bx\varphi (n))^{j+1}}%
f_{m}^{(j)}(f_{x\varphi (n)}(0))\sim \left( -1\right) ^{j}j!.
\end{eqnarray*}%
Hence the lemma follows. \hfill \rule{2mm}{3mm}\vspace{4mm}

Let%
\begin{equation*}
\mathcal{H}(n):=\left\{ 0<Z(n)\leq B\varphi (n)\right\} .
\end{equation*}

\begin{lemma}
\label{L_local}If the conditions of Theorem \ref{T_main} are valid and $%
\varphi (n)=o(n)$ as $n\rightarrow \infty ,$ then%
\begin{equation*}
\mathbf{P}\left( \mathcal{H}(n)|Z(0)=1\right) \sim \frac{\varphi (n)}{n^{2}B}%
\text{.}
\end{equation*}
\end{lemma}

\textbf{Proof}. It is known \cite{NW2006} that if the conditions of Theorem %
\ref{T_main} are valid and $k,n\rightarrow \infty $ in such a way that the
ratio $k/n$ remains bounded then
\begin{equation}
\lim_{n\rightarrow \infty }n^{2}B^{2}\left( 1+\frac{1}{Bn}\right) ^{k+1}%
\mathbf{P}\left( Z(n)=k|Z(0)=1\right) =1.  \label{Local}
\end{equation}%
Therefore,%
\begin{eqnarray*}
\mathbf{P}\left( \mathcal{H}(n)|Z(0)=1\right) &=&\sum_{1\leq k\leq B\varphi
(n)}\mathbf{P}\left( Z(n)=k|Z(0)=1\right) \\
&\sim &\frac{1}{n^{2}B^{2}}\sum_{1\leq k\leq B\varphi (n)}1\sim \frac{%
\varphi (n)}{n^{2}B}
\end{eqnarray*}%
as desired. \hfill \rule{2mm}{3mm}\vspace{4mm}

The next lemma is crucial for the proof of Theorem \ref{T_main}.

\begin{lemma}
\label{L_probab}Under the conditions of Theorem \ref{T_main} for any $x\in
\left( 0,\infty \right) $ and any $j\geq 1$%
\begin{equation*}
\lim_{n\rightarrow \infty }\mathbf{P}\left( Z(n-x\varphi (n),n\right) =j|%
\mathcal{H}(n))=\frac{x}{(j-1)!}\int_{0}^{1/x}z^{j-1}e^{-z}dz.
\end{equation*}
\end{lemma}

\textbf{Proof}. Clearly, for any $j\geq 1$%
\begin{eqnarray}
\mathbf{P}\left( Z(m,n\right) =j)&=&\sum_{k=j}^{\infty }\mathbf{P}\left(
Z(m)=k;Z(m,n\right) =j)  \notag \\
&=&\sum_{k=j}^{\infty }\mathbf{P}\left( Z(m)=k\right)
C_{k}^{j}f_{n-m}^{k-j}(0)\left( 1-f_{n-m}(0)\right)^{j}  \notag \\
&=&\frac{\left( 1-f_{n-m}(0)\right)^{j}}{j!}f_{m}^{(j)}(f_{n-m}(0)).
\label{ExpliciteZmn}
\end{eqnarray}
This representation, (\ref{SurvivalProbab}) and Lemma \ref{L_derivative} give%
\begin{eqnarray}
\mathbf{P}\left( Z(n-x\varphi (n),n\right) =j)&=&\frac{\left( 1-f_{x\varphi
(n)}(0)\right)^{j}}{j!}f_{n-x\varphi (n)}^{(j)}\left( f_{x\varphi
(n)}(0)\right)  \notag \\
&\sim &\frac{1}{j!\left( xB\varphi (n)\right)^{j}}\frac{j!\left( xB\varphi
(n)\right)^{j+1}}{B^{2}n^{2}}\sim \frac{x\varphi (n)}{Bn^{2}}.
\label{ExplicitZ}
\end{eqnarray}
Let now $Z_{1}^{\ast }(m),\ldots ,Z_{j}^{\ast }(m)$ be i.i.d. random
variables distributed as $\left\{ Z(m)|Z(m)>0\right\} ,$ and let $\eta
_{1},\ldots ,\eta _{j}$ be i.i.d. random variables having exponential
distribution with parameter 1. It is not difficult to understand, using (\ref%
{Yaglom}) that
\begin{eqnarray}
&&\lim_{n\rightarrow \infty }\mathbf{P}\left( \mathcal{H}(n)|Z(n-x\varphi
(n),n)=j\right)  \notag \\
&=&\lim_{n\rightarrow \infty }\mathbf{P}\left( Z_{1}^{\ast }(x\varphi
(n))+\cdots +Z_{j}^{\ast }(x\varphi (n))\leq B\varphi (n)\right)  \notag \\
&=&\lim_{n\rightarrow \infty }\mathbf{P}\left( \frac{Z_{1}^{\ast }(x\varphi
(n))}{Bx\varphi (n)}+\cdots +\frac{Z_{j}^{\ast }(x\varphi (n))}{Bx\varphi (n)%
}\leq \frac{1}{x}\right)  \notag \\
&=&\mathbf{P}\left( \eta _{1}+\cdots +\eta _{j}\leq \frac{1}{x}\right) =%
\frac{1}{(j-1)!}\int_{0}^{1/x}z^{j-1}e^{-z}dz.  \label{Ba}
\end{eqnarray}%
Combining this result with Lemma \ref{L_local} and (\ref{ExplicitZ}) we see
that%
\begin{eqnarray*}
&&\mathbf{P}\left( Z(n-x\varphi (n),n)=j|\mathcal{H}(n)\right) \\
&=&\frac{\mathbf{P}\left( Z(n-x\varphi (n),n\right) =j)\mathbf{P}\left(
\mathcal{H}(n)|Z(n-x\varphi (n),n)=j\right) }{\mathbf{P}\left( \mathcal{H}%
(n)\right) } \\
&\sim &\frac{x\varphi (n)}{Bn^{2}}\frac{n^{2}B}{\varphi (n)}\frac{1}{(j-1)!}%
\int_{0}^{1/x}z^{j-1}e^{-z}dz=\frac{x}{(j-1)!}\int_{0}^{1/x}z^{j-1}e^{-z}dz.
\end{eqnarray*}%
\hfill \rule{2mm}{3mm}\vspace{4mm}

\textbf{Proof of Theorem \ref{T_main}}. By the dominated convergence theorem
we have
\begin{eqnarray*}
\lim_{n\rightarrow \infty }\mathbf{E}\left[ s^{Z(n-x\varphi (n),n)}|\mathcal{%
H}(n)\right] &=&\sum_{j=1}^{\infty }\lim_{n\rightarrow \infty }\mathbf{P}%
\left( Z(n-x\varphi (n),n\right) =j|\mathcal{H}(n))s^{j} \\
&=&\sum_{j=1}^{\infty }\frac{x}{(j-1)!}\int_{0}^{1/x}s^{j}z^{j-1}e^{-z}dz \\
&=&xs\int_{0}^{1/x}e^{(s-1)z}dz=\frac{xs}{1-s}\left( 1-e^{-(1-s)/x}\right) .
\end{eqnarray*}%
Theorem \ref{T_main} is proved. \hfill \rule{2mm}{3mm}\vspace{4mm}

\textbf{Proof of Corollary \ref{Cor1}}. Since
\begin{equation*}
\mathbf{P}\left( d(n)\leq x\varphi (n)|\mathcal{H}(n)\right) =\mathbf{P}%
\left( Z(n-x\varphi (n),n)=1|\mathcal{H}(n)\right) ,
\end{equation*}%
the desired statement follows from Lemma \ref{L_probab} with $j=1.$ \hfill%
\rule{2mm}{3mm}\vspace{4mm}

\section{Proof of Theorem \protect\ref{T_simple}}

Similarly to (\ref{Ba}) we have
\begin{eqnarray}
&&\lim_{n\rightarrow \infty }\mathbf{P}\left( 0<Z(n)\leq aBn|Z(nt,n)=j\right)
\notag \\
&=&\lim_{n\rightarrow \infty }\mathbf{P}\left( Z_{1}^{\ast }(n(1-t))+\cdots
+Z_{j}^{\ast }(n(1-t))\leq aBn\right)  \notag \\
&=&\lim_{n\rightarrow \infty }\mathbf{P}\left( \frac{Z_{1}^{\ast }(n(1-t))}{%
Bn(1-t)}+\cdots +\frac{Z_{j}^{\ast }(n(1-t))}{Bn(1-t)}\leq \frac{a}{1-t}%
\right)  \notag \\
&=&\mathbf{P}\left( \eta _{1}+\cdots +\eta _{j}\leq \frac{a}{1-t}\right) =%
\frac{1}{(j-1)!}\int_{0}^{a/(1-t)}z^{j-1}e^{-z}dz.  \label{Bo2}
\end{eqnarray}%
Besides,%
\begin{equation*}
\lim_{n\rightarrow \infty }\mathbf{P}\left( Z(n)\leq aBn|Z(n)>0\right)
=1-e^{-a}
\end{equation*}%
and, by (\ref{ExpliciteZmn}), (\ref{SurvivalProbab}) and Lemma \ref%
{L_Deriv_n}
\begin{eqnarray*}
\mathbf{P}\left( Z(nt,n\right) =j) &=&\frac{\left( 1-f_{n(1-t)}(0)\right)
^{j}}{j!}f_{nt}^{(j)}(f_{n(1-t)}(0)) \\
&\sim &\left( \frac{1}{Bn\left( 1-t\right) }\right) ^{j}\frac{1}{j!}\frac{%
j!\left( \frac{1-t}{t}\right) ^{2}\left( Bn\left( 1-t\right) \right) ^{j-1}}{%
\left( \frac{1-t}{t}+1\right) ^{j+1}} \\
&\sim &\frac{1-t}{Bn}t^{j-1}.
\end{eqnarray*}%
Therefore,%
\begin{eqnarray}
&&\lim_{n\rightarrow \infty }\mathbf{P}\left( Z(nt,n)=j|0<Z(n)\leq aBn\right)
\notag \\
&=&\lim_{n\rightarrow \infty }\frac{\mathbf{P}\left( 0<Z(n)\leq
aBn|Z(nt,n)=j\right) \mathbf{P}\left( Z(nt,n)=j\right) }{\mathbf{P}\left(
Z(n)\leq aBn|Z(n)>0\right) \mathbf{P}\left( Z(n)>0\right) }  \notag \\
&=&\frac{1-t}{1-e^{-a}}\frac{1}{\left( j-1\right) !}\int_{0}^{a/(1-t)}\left(
tz\right) ^{j-1}e^{-z}dz.  \label{Distant_simple}
\end{eqnarray}%
As a result we get%
\begin{eqnarray*}
\lim_{n\rightarrow \infty }\mathbf{E}\left[ s^{Z(nt,n)}|0<Z(n)\leq aBn\right]
&=&\frac{(1-t)s}{1-e^{-a}}\sum_{j=1}^{\infty }\int_{0}^{a/(1-t)}\frac{\left(
stz\right) ^{j-1}}{\left( j-1\right) !}e^{-z}dz \\
&=&\frac{(1-t)s}{1-e^{-a}}\int_{0}^{a/(1-t)}e^{(st-1)z}dz \\
&=&\frac{(1-t)s}{(1-e^{-a})(1-ts)}\left( 1-e^{-(1-ts)a/(1-t)}\right) .
\end{eqnarray*}%
Theorem \ref{T_simple} is proved. \hfill \rule{2mm}{3mm}\vspace{4mm}

\textbf{Proof of Corollary \ref{Cor2}}. Since
\begin{equation*}
\mathbf{P}\left( d(n)\leq tn|0<Z(n)\leq aBn\right) =\mathbf{P}\left(
Z(n(1-t),n)=1|0<Z(n)\leq aBn\right) ,
\end{equation*}%
the desired statement follows from (\ref{Distant_simple}) with $j=1$ and $%
1-t $ for $t$. \hfill \rule{2mm}{3mm}


\vspace{4mm}

\begin{thebibliography}{99}
\bibitem{Athr2012} Athreya, K. B. (2012) Coalescence in the recent past in
rapidly growing populations. \textit{Stochastic Processes and their
Applications,} \textbf{122}, 3757--3766.

\bibitem{Ath2012b} Athreya, K. B. (2012) Coalescence in critical and
subcritical Galton-Watson branching processes. \textit{Journal of Applied
Probability,} \textbf{49}, 627--638.

\bibitem{AN72} Athreya, K. B. and Ney, P. E. (1972) \textit{Branching
processes.} Springer--Verlag, Berlin-Heidelberg--New York .

\bibitem{Dur78} Durrett, R. (1978) The genealogy of critical branching
processes. \textit{Stochastic Processes and their Applications,} \textbf{8},
101--116.

\bibitem{FP} Fleischmann, K., Prehn, U. (1974) Ein Grenzfersatz f\"{u}r
sub\-kri\-tische Ver\-zwei\-gungs\-pro\-zes\-se mit eindlich vielen Typen
von Teilchen. \textit{Math. Nachr.}, \textbf{64}, 233--241.

\bibitem{FZ} Fleischmann, K., Siegmund-Schultze, R. (1977) The structure of
reduced critical Galton-Watson processes. \textit{Math. Nachr.}, \textbf{79}%
, 233--241.

\bibitem{HJR17} Harris, S. C., Johnston, S. G. G., and Roberts, M. I. (2017)
The coalescent structure of continuous-time Galton-Watson trees.
https://arxiv.org/pdf/1703.00299.pdf

\bibitem{Sam17} Johnston, S. G. G. (2017) Coalescence in supercritical and
subcritical continuous-time Galton-Watson trees.
https://arxiv.org/pdf/1709.008500v1.pdf

\bibitem{Lam2003} Lambert, A. (2003) Coalescence times for the branching
process. \textit{Advances in Applied Probability,} \textbf{35}, 1071--1089.

\bibitem{Lam2016} Lambert, A. (2016) Probabilistic models for the subtrees
of life. https://arxiv.org/abs/1603.03705

\bibitem{Le2014} Le, V. (2014) Coalescence times for the Bienaym\'{e}%
-Galton-Watson process. \textit{Journal of Applied Probability,} \textbf{51}%
, 209--218.

\bibitem{NW2006} Nagaev, S. V. and Vakhtel, V. I. (2006) On the local limit
theorem for a critical Galton--Watson process. \textit{Theory Probab. Appl.}%
, \textbf{50}, 400--419.

\bibitem{Con95} O'Connell, N. (1995) The genealogy of branching processes
and the age of our most recent common ancestor. \textit{Advances in Applied
Probability}, \textbf{27}, 418--442.

\bibitem{Sen67} Seneta, E. (1967) The Galton-Watson process with mean one.
\textit{Journal of Applied Probability, }\textbf{4}, 489--495.

\bibitem{Sev74} Sewast'yanov, B. A. (1974) Verzweigungsprozesse.
Mathematische Lehrbucher und Monographien. II. \textit{Abteilung:
Mathematische Monographien, Band 34. Akademie-Verlag, Berlin,} xi+326 .

\bibitem{VD2008} Vatutin, V. A., Dyakonova, E. E. (2008) Limit theorems for
reduced processes in random environment. \textit{Theory Probab. Appl.},
\textbf{52}, 277--302.

\bibitem{VD2015} Vatutin, V. A. and D'yakonova, E. E. (2015) Decomposable
branching processes with a fixed extinction moment. \textit{Proc. Steklov
Inst. Math.}, \textbf{290}, 103--124.

\bibitem{Zub} Zubkov, A. M. (1975) Limit distributions of the distance to
the nearest common ancestor. \textit{Theory Probab. Appl.}, \textbf{20},
602--612.
\end{thebibliography}
\end{document}